%% Alexander Tumanov
%%
%% TESTING ANALYTICITY ON CIRCLES
%%
%% Plain TeX

%\input Srctex.sty

\magnification=\magstep1
\baselineskip=15pt
\parskip=4pt

\def\a{\alpha}
\def\b{\beta}
\def\C{{\bf C}}
\def\R{{\bf R}}
\def\A{{\cal A}}

\def\jour{\sl}
\def\vol{\bf}
\def\bysame{\vrule height0.4pt depth0pt width1truein \ , }

\topglue 1truein
\centerline{\bf TESTING ANALYTICITY ON CIRCLES}
\bigskip
\centerline{A. Tumanov}
\bigskip

Let $\Omega$ be a domain in complex plane
and let $C_t$ be a continuous one parameter family
of Jordan curves such that $\cup C_t = \bar\Omega$.
Let $f$ be a continuous function in
$\bar\Omega$ such that the restrictions
$f|_{C_t}$ extend holomorphically inside $C_t$.
When does this imply that $f$ is holomorphic in $\Omega$?
This question has been known for a long time
and is related to inverse problems for PDE
and integral geometry, see [Z,E1].

Globevnik [G] answered the question in the affirmative for
rotation invariant families. Agranovsky and Globevnik [AG]
resolved the question for families of circles when $f$ is a
rational or real-analytic functions of two real variables.
Ehrenpreis [E2] independently found solution for the circles of
constant radii with centers on a line and real-analytic functions.
In [T] we obtain the affirmative answer for the family of
radius one circles $C_t$ with centers at $t\in\R$, where
$|t|<2+\epsilon$, $\epsilon>0$ and for two more special families.
Ehrenpreis [E1] has independently obtained the main result of [T].

In this paper we extend the result of [T] to fairly
general families of circles with variable radii.
In particular, the conclusion holds for radius one
circles with centers at $t\in\R$, $|t|<1+\epsilon$,
improving the result of [T].
However for curves other than the circles the
question is largely open.

Despite the one variable nature of the problem,
our method involves analysis of several complex variables,
in particular the extendibility of CR functions.
Surprisingly, the proof is simpler than in our
previous paper. The main tool comes form the original
proof of the classical H. Lewy [L] extension theorem.

The author thanks Mark Agranovsky, Leon Ehrenpreis,
and Josip Globevnik for useful discussions.
\bigskip
\centerline{*\quad*\quad*}
\medskip

Let $\{C_t : \a\le t \le\b\}$, be a continuous one parameter
family of circles in complex plane $\C$ with centers at $c(t)\in\C$
and radii $r(t)>0$. Let $D_t$ denote the disc bounded by $C_t$.
Suppose the following hold.

\item{(a)} $\bar D_\a\cap\bar D_\b=\emptyset$, that is
$|c(\a)-c(\b)|>r(\a)+r(\b)$.

\item{(b)} The functions $c(t)$ and $r(t)$ are piecewise $C^3$
smooth. The curve $t\mapsto c(t)$ is injective
and regular, that is $c'(t)\ne 0$.

\item{(c)} No circle $C_t$ is contained in the closed disc
$\bar D_s$ for $t\ne s$, that is
$|c(t)-c(s)|>|r(t)-r(s)|$.

\item{(d)} $|c'(t)|>|r'(t)|$ whenever defined.
\medskip

{\bf Theorem.}
Let the family $\{C_t : \a\le t \le\b\}$ satisfy (a)-(d).
Let $\Omega=\bigcup D_t$.
Let $f:\bar\Omega\to\C$ be a continuous function.
Suppose for every $\a\le t\le\b$ the restriction
$f|_{C_t}$ extends holomorphically to $D_t$.
Then $f$ is holomorphic in $\Omega$.
\medskip

{\it Remarks.} The condition (a) is crucial, the others being added
for simplicity and convenience of the proof.
The smoothness in excess of $C^1$ is used only to
deal with triple intersections of the circles.
The condition that the circles can't lie inside one another
is natural because otherwise the values of $f$ on them
are unrelated. We assume the slightly stronger
property (c) that they can't even touch. The condition (d)
is the infinitesimal version of (c). In fact (c) implies
(d) with possible equality, but for simplicity we assume
the strict inequality.

Turning to the proof, we define
$$
\Sigma=\{(z,w)\in\C^2: w=\bar z \}
$$$$
X_t=\{(z,w)\in\C^2: (z-c(t))(w-\bar c(t))=r(t)^2,
|z-c(t)|\le r(t) \}.
$$
The complex curve $X_t$ can be considered (a part of)
the complexification of $C_t$. Indeed, let
$\tilde C_t=\{(z,w)\in\C^2: z\in C_t, w=\bar z\}$.
Then $\partial X_t=\tilde C_t$ because
$(z,w)\in X_t$ and $|z-c(t)|= r(t)$ imply
$w-\bar c(t)={r(t)^2\over z-c(t)}=\bar z-\bar c(t)$,
whence $w=\bar z$.
\medskip

{\bf Lemma 1.
\it
The condition (c) implies
$X_t\cap X_s=\tilde C_t\cap\tilde C_s\subset\Sigma$ for $t\ne s$.

Proof.}
By (c) we have the following possibilities.

Case 1: $\bar D_t\cap\bar D_s=\emptyset$.
Then $X_t\cap X_s=\emptyset$ because
$(z,w)\in X_t\cap X_s$ implies $z\in\bar D_t\cap\bar D_s$.

Case 2: $C_t\cap C_s\ne\emptyset$.
By eliminating $w$ from the equations of $X_t$ and $X_s$,
we get a quadratic equation in $z$. Hence $X_t\cap X_s$
contains no more points than $\tilde C_t\cap\tilde C_s$.
The lemma is proved.
\medskip

Define
$M=\bigcup X_t$. Then by Lemma 1, $M\setminus\Sigma$ is a
piecewise smooth Levi-flat hypersurface in $\C^2$.
Let $f_t$ denote the holomorphic extension of $f$
to $D_t$. For $(z,w)\in X_t$ we define
$F(z,w)=f_t(z)$. Then $F$ is a continuous CR function
on $M$ because $F$ is holomorphic on the fibers $X_t$.

We plan to prove that $F$ actually is independent of $w$.
That would mean that all the extensions $f_t(z)$
match at $z$, and $f$ is holomorphic.

Let
$\Omega'=\Omega\setminus(\bar D_\a\cup\bar D_\b)$.
For $z\in\Omega'$ put
$I_z=\{\a\le t\le\b: |z-c(t)|\le r(t)\}\subset\R$.
Since $z\in\Omega'$, then $\a$ and $\b$ are not in $I_z$.
Let $\bar\C$ be the Riemann sphere $\C\cup\{\infty\}$.

Put $\Gamma_z=\{w\in\bar\C: (z,w)\in M \}$.
If $\Gamma_z$ is unbounded, then we will
consider $\infty\in\Gamma_z$.
It will suffice to consider the case in which $I_z$
consists of finitely many closed disjoint intervals
$I_z=I_1\cup ... \cup I_k$, where $k$ depends on $z$.
Then $\Gamma_z$ is a parametrized curve $I_z\to \bar\C$,
$t\mapsto w(t)=\bar c(t)+{r(t)^2\over z-c(t)}$.

If $t$ is an end point of one of the intervals $I_j$,
then $|z-c(t)|=r(t)$ and as we noticed before,
$w(t)=\bar z$. By Lemma 1 the mapping
$t\mapsto w(t)$ is injective except that all the end
points of the intervals $I_j$ are mapped to the same
point $\bar z$.

Hence the curve $\Gamma_z$ consists of $k$ simple
closed loops
$\Gamma_z= \Gamma_1+ ... +\Gamma_k$
corresponding to the partition
$I_z=I_1\cup ... \cup I_k$.

Let $C=\{c(t):\a\le t\le\b\}$.
If $z\in C$, then exactly one loop of the above
passes through the point at infinity.

We first explain the idea of the proof in the
case when
the circles have no triple intersections.
Then the curve $\Gamma_z$ consists of just one
loop. By the argument of the classical work
of H. Lewy [L], the CR function $F$
holomorphically extends inside $\Gamma_z$
for $z\in\Omega'\setminus C$.
When $z$ crosses the centers line $C$,
the interior and exterior of the loop
$\Gamma_z$ interchange. This implies that
$F(z,w)$ extends to the whole Riemann sphere
for fixed $z$, hence $F(z,w)$ is independent of
$z$, and $f$ is holomorphic.

In the general case
we examine the behavior of $\Gamma_z$
as a function of $z$.
Consider the natural map
$Z:(\a,\b)\times\R\to\C$,
$Z(t,\theta)=c(t)+r(t)e^{i\theta}$.
By the implicit function theorem, if $z_0$
is not a critical value of $Z$, then for
$z$ close to $z_0$,
$\Gamma_z$ varies continuously with $z$,
in particular
the number of loops in $\Gamma_z$ is constant.

Let $P$ denote the set of all critical values of $Z$.
They are found by solving the equation
${\partial Z\over\partial t} =
\lambda {\partial Z\over\partial\theta}$,
$\lambda\in\R$.
This equation means that the velocity of the point
of $C_t$ is tangential to $C_t$, so we call such
point {\it sliding} points.
Straightforward calculations show that
every circle has exactly two sliding points
$$
p(t)=c(t)-
{r'(t)\pm i\sqrt{|c'(t)|^2-r'(t)^2}\over \bar c'(t)}.
$$
Hence $P$ consists of finitely many $C^2$-smooth
curves possibly with singular points,
in which $p'(t)=0$.
Regular pieces of $P$ are the osculating curves
for the family $\{C_t\}$, that is $C_t$ is tangent
to $P$ at $p(t)$.

Let $P'$ be the set of regular points of $P$.
Without loss of generality $P'$ is ``simple",
that is if $p(t)$ and $p(s)$ are regular points
for $t\ne s$, then $p(t)\ne p(s)$.
Indeed, by (c) the circles $C_t$ and $C_s$ cannot
touch in interior fashion. Exterior tangency
can be eliminated by shrinking the interval
$(\a,\b)$ and passing to a subfamily.

We describe the qualitative behavior of $\Gamma_z$
as $z$ crosses the critical set $P$ at $p(t)$.
Let $p(t)\in P'$.
Let $\rho(t)$ denote the radius of curvature of $P$
at $p(t)$.
We distinguish between the interior and exterior
tangency of $P$ and $C_t$ at $p(t)$ if
$\rho(t)<\infty$.
Only the following cases may occur.

\item{(1)} Interior tangency and $\rho(t)>r(t)$ or
exterior tangency.

\item{(2)} Interior tangency and $\rho(t)<r(t)$.

\item{(3)} Interior tangency and $\rho(t)=r(t)$.

The last case (3) actually can't occur because it
implies $c'(t)=0$ in violation of (b).
In case (1), if z crosses to the side of $P$ where
$C_t$ approaches $p(t)$, then a new small loop
of $\Gamma_z$ is created.
In case (2), one loop of $\Gamma_z$ splits into
two loops.

As a tool for proving that $F(z,w)$ is constant on
$\Gamma_z$, we use the following simple
\medskip

{\bf Lemma 2.
\it
Let $G\subset\bar\C$ be a closed piecewise smooth
curve and let $f\in L^\infty(G)$.
Suppose
$\int_G f(\zeta)(\zeta-w)^{-2}\,d\zeta=0$
for all $w\in\C\setminus G$.
Then $f$ is constant on every segment of $G$
with no points of self-intersection of $G$.

Proof.}
Consider the Cauchy type integral
$F(w)={1\over 2\pi i}
\int_G {f(\zeta)\,d\zeta\over \zeta-w}$,
$w\in\C\setminus G$.
The hypotheses imply that $F'(w)\equiv 0$,
whence $F$ is locally constant.
The lemma now follows by
the Plemelj-Sokhotsky jump formula.
\medskip

In view of the last lemma, we put
$$
\Phi(z,w)=\int_{G_z}F(z,\zeta)
(\zeta-w)^{-2}d\zeta,
$$
where $G_z$ is a curve that consists
of some of the loops of $\Gamma_z$ and depends
continuously on $z$ in an open set
$U\subset\Omega'$, and $w\in\C\setminus G_z$.
The second power in the denominator is convenient
to avoid convergence problems if $G_z$
passes though the infinity.
\medskip

{\bf Lemma 3.
\it
$\Phi(z,w)$ is holomorphic in
$z\in\Omega'\setminus(P\cup C)$,
$w\in\C\setminus G_z$.

Proof.}
Obviously, $\Phi$ is holomorphic in $w$.
To show that $\Phi$ is holomorphic in $z$,
we follow the original proof of the
H. Lewy [L] extension theorem.
Put $\Phi(z)=\Phi(z,w)$,
$H(z,\zeta)=F(z,\zeta)(\zeta-w)^{-2}$, then
$\Phi(z)=\int_{G_z}H(z,\zeta)\,d\zeta$.
By the Morera theorem, it suffices to show
that $\int_\gamma\Phi(z)\,dz=0$
for every small loop $\gamma$ in $U$.

Without loss of generality, $G_z$ consists
of just one simple loop. Consider the
``torus''
$T=\{(z,\zeta)\in M:
z\in\gamma, \zeta\in G_z\}$.
Then $T$ bounds a ``solid torus''
$S\subset M$ obtained by filling the loop
$\gamma$.
Then by Stokes' formula
$$
\int_\gamma\Phi(z)\,dz=
\int_T H(z,\zeta)\,d\zeta\wedge dz=
\int_S dH(z,\zeta)\wedge d\zeta\wedge dz=0
$$
because $H$ is a CR function.
The lemma is proved.
\medskip

Let $L_0$ be a straight line separating the
circles $C_\a$ and $C_\b$.
Let $L$ be a small perturbation of $L_0$
such that

\item{(i)} $L \cap (P\setminus P')=\emptyset$;
\item{(ii)} $L\cap (P\cup C)$ is finite;
\item{(iii)} $L$ has only transverse
intersections with $P$ and $C$.

Such a curve $L$ exists because
$P$ and $C$ are smooth and the singular set
$P\setminus P'$ has length zero.

Consider the set $\A$ of all loops of
$\Gamma_z$ for all $z\in L$.
For $G_1,G_2\in\A$,
$G_1\subset\Gamma_{z_1}$,
$G_2\subset\Gamma_{z_2}$,
$z_1,z_2\in L$,
we put $G_1\sim G_2$ if
$G_1$ deforms into $G_2$ as $z$ runs from
$z_1$ to $z_2$ in $L$.
If a single loop $G\in\A$ splits into two loops
$G_1$ and $G_2$ as $z\in L$ crosses the critical
set $P$, then we also put $G_1\sim G_2$.
We define an equivalence relation $\sim$ on $\A$
using the above two basic equivalences.

Since the centers line $C$ comes from one side
of $L$ to the other and all intersections $L\cap C$
are transverse, then $L\cap C$ consists of {\it odd}
number of points.
Since no more than one loop of $\Gamma_z$ passes
through the infinity, then there exists an
equivalence class $A\subset\A$ that includes
{\it odd} number of infinite loops.

Define $G_z$ as a curve that consists of all
loops of $\Gamma_z$ that are in $A$.

Let $U_0$ be a small neighborhood of $L$.
Then $U_0$ is an infinite thin band separated into
sub-bands by segments of the singular set $P$ and
centers line $C$. We continuously extend $G_z$
from $z\in L$ to $z\in U_0$. We restrict to
the points $z$ for which $G_z$ is nontrivial,
that is we put
$U=\{z\in U_0:G_z\ne\emptyset,
G_z\ne\{\bar z\}\}$.

By Lemma 3 the function $\Phi(z,w)$ defined above
is holomorphic with respect to both $z$ and $w$
for $z\in U\setminus(P\cup C)$,
$w\in\bar\C\setminus G_z$.
Since $G_z$ varies continuously and the integrand
in the formula for $\Phi$ decays as $\zeta^{-2}$
at infinity, then $\Phi$ is continuous
on $V=\{(z,w): z\in U, w\in\bar\C\setminus G_z\}$.
Hence $\Phi$ is holomorphic on $V$.

Let $G$ be an oriented piecewise smooth curve
in $\bar\C$. Define $G^+$ and $G^-$ as the sets of
all points in $\bar\C$ of index respectively
1 and 0 with respect to $G$. We call $G$
{\it quasi-simple} if
$G^+\cup G^-\cup G=\bar\C$.

By the nature of the equivalence relation $\sim$,
all the curves $G_z$ are quasi-simple
(which might not be the case for $\Gamma_z$), so
$V=V^+\cup V^-$, where
$V^\pm=\{(z,w): z\in U, w\in G_z^\pm \}$.
Also, it follows that the sets $V^\pm$
are connected.

The set $U$ is a thin band with two short edges
$K^\pm\subset P$. If $z\in U$ is close to $K^\pm$,
then $G_z$ consists of a single small loop
contracting into $\bar z$ as $z$ approaches
$K^\pm$.

Here is the key point of the proof.
As $z\in U$ runs from $K^+$ to $K^-$,
it crosses the centers line $C$ odd number of
times so $G_z$ crosses through the $\infty$
odd number of times.
If say $\infty\in G_z^+$
for $z\in U$ close to $K^+$,
then $\infty\in G_z^-$
for $z\in U$ close to $K^-$.
For definiteness choose the signs in
the notation $K^\pm$ so this is the case.

As $z$ approaches $K^\pm$, the loop $G_z$ contracts
into a point, hence
$\Phi(z,w)\to 0$ as $z\to z_0\in K^\pm$.
By the uniqueness theorem,
$\Phi\equiv0$ in the connected set $V^\pm$.
By Lemma 2 the function $F$ is constant on every
loop of $G_z$.
This implies that $f$ is holomorphic in $U$
and this property propagates over the whole
set $\Omega$.
The theorem is now proved.

\beginsection References

\frenchspacing

\item{[AG]}
M. Agranovsky and J. Globevnik,
Analyticity on circles for rational and real analytic
functions of two real variables,
{\jour J. D'Analyse Math. \vol 91},
(2003),31--65.

\item{[E1]}
L. Ehrenpreis,
{\jour The Universality of the Radon Transform},
Oxford Univ. Press (2003).

\item{[E2]}
\bysame
Three problems at Mount Holyoke,
{\jour Contemp. Math., \vol 278},
Providence, RI, 2001.

\item{[G]}
J. Globevnik,
Testing analyticity on rotation invariant families of curves,
{\jour Trans. Amer. Math. Soc. \vol 306}
(1988), 401--410.

\item{[L]}
H. Lewy,
On the local character of the solutions of an atypical
linear differential equation in three variables and a
related theorem for regular functions of two complex
variables,
{\jour Ann. Math. \vol 64}
(1956), 514--522.

\item{[T]}
A. Tumanov,
A Morera type theorem in the strip,
{\sl Math. Res. Let. \vol 11} (2004), 23--29.

\item{[Z]}
L. Zalcman,
Analyticity and the Pompeiu problem,
{\jour Arch. Rat. Mech. Anal. \vol 47}
(1972), 237--254.

\bigskip\noindent
Alexander Tumanov, Department of Mathematics, University of Illinois,
Urbana, IL 61801. E-mail: tumanov@math.uiuc.edu
\bye